\documentclass[reqno]{amsart}
\usepackage{amssymb}
\usepackage{hyperref}

\numberwithin{equation}{section}
\newtheorem{theorem}{Theorem}[section]
\newtheorem{lemma}[theorem]{Lemma}

\newtheorem{corollary}[theorem]{Corollary}

\begin{document}
\title{ Some spectral properties for generalized derivations }

\author[ M. Amouch and F. Lombarkia ]
{ M. Amouch  and F. Lombarkia }  
\address{Mohamed AMOUCH \newline
Department of Mathematics 
University Chouaib Doukkali,
Faculty of Sciences, Eljadida.
24000, Eljadida, Morocco.}
\email{mohamed.amouch@gmail.com}

\address{Farida Lombarkia \newline
 Department of Mathematics, Faculty of Science, University of Batna,
05000, Batna, Algeria.}
\email{lombarkiafarida@yahoo.fr}



\subjclass[2000]{47A10, 47A53, 47B47}
\keywords{left polaroid; elementary operator; finitely left polaroid}

\begin{abstract}
Given Banach spaces $\mathcal{X}$ and $\mathcal{Y}$
and Banach space operators $A\in L(\mathcal{X})$ and $B\in L(\mathcal{Y}).$
The generalized derivation $\delta_{A,B} \in L(L(\mathcal{Y},\mathcal{X}))$
is defined by $\delta_{A,B}(X)=(L_{A}-R_{B})(X)=AX-XB$.
This paper is concerned with the problem of the transferring
the left polaroid property, from  operators $A$ and $B^{*}$
to the generalized derivation $\delta_{A,B}$.
As a consequence, we give necessary and sufficient conditions for $\delta_{A,B}$ to satisfy
generalized a-Browder's theorem and generalized a-Weyl's theorem. As application, we extend some recent
results concerning Weyl type theorems.
\end{abstract}
\maketitle

\section{Introduction}
Given Banach spaces $\mathcal{X}$ and $\mathcal{Y}$
and Banach space operators $A\in L(\mathcal{X})$ and $B\in L(\mathcal{Y}).$
Let $L_{A}\in L(L(\mathcal{X}))$ and $R_{B}\in L(L(\mathcal{Y}))$ be the left and the right
multiplication operators, respectively, and denote by $\delta_{A,B}\in L(L(\mathcal{Y},\mathcal{X}))$
the generalized derivation $\delta_{A,B}(X)=(L_{A}-R_{B})(X)=AX-XB$. The problem of transferring spectral properties from
$A$ and $B$ to $L_A,$ $R_B$, $L_A R_B$ and $\delta_{A,B}$ was studied by numerous mathematicians,
see \cite{Bao, BD, BA, CH, D2, D9, H2, LB, LBA} and the references therein.
The main objective of this paper is to study the
problem of transferring the left polaroid property and its strong version finitely left polaroid property,
from $A$ and $B^{*}$ to $\delta_{A,B}.$
After section 2 where several basic definitions and facts will be recalled,
we will prove that if $A$ is left polaroid and satisfies property $(\mathcal{P}_l)$ and $B$ is right polaroid and
satisfy property $(\mathcal{P}_r),$ then $\delta_{A,B}$ is left polaroid.
Also, we prove that $A$ is finitely left polaroid and $B$ is finitely right polaroid
if and only if $\delta_{A,B}$ is finitely left polaroid.
In the fourth section, we give necessary and sufficient conditions for $\delta_{A,B}$ to satisfy generalized a-Weyl's theorem.
In the last section we apply results obtained previously. If $\mathcal{X}=H$ and $\mathcal{Y}=K$ are Hilbert spaces,
we prove that if $A\in L(H)$ and $B\in L(K)$ are
completely totally hereditarily normaloid operators, then
$f(\delta_{A,B})$ satisfies generalized a-Weyl's theorem, for every analytic function $f$ defined on a neighborhood of
$\sigma(\delta_{A,B})$ which is non constant on each of the components of its domain. This generalize results obtained
in \cite{BA, CH, D2, Du8, LB, LBA}.

\section{Notation and terminology}
Unless otherwise stated, from now on $\mathcal{X}$ (similarly, $\mathcal{Y}$)
shall denote a complex Banach space and $L(\mathcal{X})$ (similarly, $L(\mathcal{Y})$)
the algebra of all bounded linear maps defined on and with values in $\mathcal{X}$
(respectively, $\mathcal{Y}$).
Given $T\in L(\mathcal{X}),$ $N(T)$ and $R(T)$ will stand for the null space and the range of $T$
respectively. Recall that $T\in L(\mathcal{X})$ is said to be bounded below, if $N(T)=\{0\}$
and $R(T)$ is closed. Denote the approximate point spectrum of $T$ by
\begin{equation*}
\sigma_{a}(T)=\{\lambda\in\mathbb{C}: T-\lambda I\,\,\,\text{is not bounded below}\}.
\end{equation*}
Let \begin{equation*}\sigma_{s}(T)=\{\lambda\in\mathbb{C}: T-\lambda I\,\,\,\text{is not surjective}\}\end{equation*}
denote the surjective spectrum of $T.$
In addition, $\mathcal{X}^*$ will denote the dual space of $\mathcal{X},$ and if
$T\in \mathcal{X},$ then $T^* \in L(\mathcal{X}^{*})$ will stand for the adjoint map
of $T.$ Clearly, $\sigma_{a}(T^*)=\sigma_{s}(T)$ and
$\sigma_{a}(T) \cup \sigma_{s}(T)=\sigma(T),$ the spectrum of $T.$ Recall that the ascent $asc(T)$ of an operator $T$ is defined by $asc(T)=\inf\{n\in\mathbb{N}:N(T^{n})=N(T^{n+1})\}$ and the descent
$dsc(T)=\inf\{n\in\mathbb{N}:R(T^{n})=R(T^{n+1})\}$, with $\inf\emptyset=\infty.$ It is well known that if $asc(T)$ and $dsc(T)$ are both finite, then they are equal.\\
\indent A complex number
$\lambda\in\sigma_{a}(T)$  (respectively, $\lambda\in\sigma_{s}(T)$) is left pole (respectively, right pole) of order $d$  of $T\in L(\mathcal{X})$ if $asc(T-\lambda I)=d<\infty$ and $R((T-\lambda I)^{d+1})$ is closed (respectively, $dsc(T-\lambda I)=d<\infty$ and $R((T-\lambda I)^{d})$ is closed).
We say that $T$ is left polar (respectively, right polar) of order d at a point $\lambda\in \sigma_a(T)$
(respectively, $\lambda \in \sigma_s(T)$) if $\lambda$ is a left pole of $T$ (respectively, right pole of $T$) of order d.
Now, $T$ is left polaroid (respectively, right polaroid) if $T$ is left polar (respectively, right polar ) at
every $\lambda \in iso\sigma_a(T)$ (respectively, $\lambda \in iso\sigma_s(T)$), where
$iso\mathcal{K}$ is the set of all isolated points of $\mathcal{K}$ for $\mathcal{K} \subseteq\mathbb{C}.$
In \cite{BD} a left polar operator $T\in L(\mathcal{X})$ of order $d(\lambda)$ at $\lambda \in \sigma_a(T),$
satisfies property $(\mathcal{P}_l)$ if the closed subspace $N((T-\lambda)^{d(\lambda)})+R(T-\lambda)$
is complemented in $\mathcal{X}$ for every $\lambda \in iso\sigma_a(T).$
Dually, a right polar operator $T\in L(\mathcal{X})$ of order $d(\lambda)$ at $\lambda \in \sigma_s(T),$
satisfies property $(\mathcal{P}_r)$ if the closed subspace
$N(T-\lambda) \cap R((T-\lambda)^{d(\lambda)})$ is complemented in $\mathcal{X}$ for every $\lambda \in iso\sigma_s(T).$
If $\mathcal{X}=H$
is a Hilbert space, then every left polar (respectively, right polar) operator $T\in L(H)$ of order $d(\lambda)$ at
$\lambda \in iso\sigma_a(T)$ (respectively, $\lambda \in iso\sigma_s(T)$)
satisfies property $(\mathcal{P}_l)$ (respectively, $(\mathcal{P}_r)).$
On the other hand, it is known that $T\in L(\mathcal{X})$ is right polaroid if and only if $T^{*}$ is left polaroid
and $T$ is polaroid if it is both left and right polaroid, whenever
$iso \sigma(T)= iso\sigma_a(T)\cup iso\sigma_s(T).$\\
\indent Recall that $T\in L(\mathcal{X})$ is said to be a Fredholm operator, if both $\alpha(T)=dimN(T)$ and
$\beta(T)=dim \mathcal{X} \setminus R(T)$ are finite dimensional, in which case its index is given by
$ind(T)=\alpha(T)-\beta(T).$ If $ R(T)$ is closed and $\alpha(T)$ is finite (respectively, $\beta(T)$ is finite),
then $T \in L(\mathcal{X})$ is said to be upper semi-Fredholm (respectively, lower semi-Fredholm)
while if $\alpha(T)$ and $ \beta(T)$ are both finite and equal, so the index is zero and $T$ is said to be Weyl operator.
These classes of opertaors generate the Fredholm spectrum, the upper semi-Fredholm spectrum, the lower semi-Fredholm spectrum
and the Weyl spectrum of $T\in L(\mathcal{X})$ which will be denoted by $\sigma_{e}(T),$
$\sigma_{SF_{+}}(T),$ $\sigma_{SF_{-}}(T)$ and $\sigma_{W}(T),$ respectively.
The Weyl essential approximate point spectrum and the Browder essential approximate point spectrum of $T\in L(\mathcal{X})$ are the sets
\begin{equation*}
\sigma_{aw}(T)=\{\lambda\in\sigma_{a}(T): \lambda \in \sigma_{SF_{+}}(T)\,\, or \,\, 0<ind(T-\lambda I)\}
\end{equation*}
and
\begin{equation*}
\sigma_{ab}(T)=\{\lambda\in\sigma_{a}(T): \lambda\in\sigma_{aw}(T)\,\,\text{
or}\,\,\ asc(T-\lambda I)=\infty\}.
\end{equation*}
It is clear that
\begin{equation*}
\sigma_{SF_{+}}(T)\subseteq\sigma_{aw}(T)\subseteq\sigma_{ab}(T)\subseteq\sigma_{a}(T).
\end{equation*}
\indent For $T\in L(\mathcal{X})$ and a nonnegative integer $n$ define $T_{n}$ to be the restriction of
$T$ to $R(T^{n})$ viewed as a map from $R(T^{n})$ into $R(T^{n})$. If for some integer $n$
the range space $R(T^{n})$ is closed and the induced operator $T_{n}\in L(R(T^{n}))$
is Fredholm, then $T$ will be said B-Fredholm.
In a similar way, if $T_{n}$ is upper semi-Fredholm (respectively, lower semi-Fredholm )
operator, then $T$ is called upper semi B-Fredholm (respectively, lower semi B-Fredholm).
In this case the index of $T$ is defined as the index of the semi-Fredholm operator $T_{n}$, see \cite{B1}.
$T\in L(\mathcal{X})$ is semi B-Fredholm if $T$ is upper semi B-Fredholm or lower semi B-Fredholm.
Let
\begin{equation*}
\Phi_{SBF}(\mathcal{X})=\{T\in L(\mathcal{X}) : T\,\,\,\text{%
is semi B-Fredholm}\,\,\},
\end{equation*}
\begin{equation*}
\Phi_{SBF_{+}^{-}}(\mathcal{X})=\{T\in\Phi_{SBF}(\mathcal{X}): T\,\,\,\text{%
is upper semi B-Fredholm with}\,\,ind(T)\leq0\}
\end{equation*}
and
\begin{equation*}
\Phi_{SBF_{-}^{+}}(\mathcal{X})=\{T\in\Phi_{SBF}(\mathcal{X}): T\,\,\,\text{%
is lower semi B-Fredholm with}\,\,ind(T) \geq 0\}.
\end{equation*}
Then the upper semi B-Weyl and lower semi B-Weyl spectrum of $T$ are the sets
\begin{equation*}
\sigma_{UBW}(T)=\{\lambda\in\sigma_{a}(T): T-\lambda I\not\in\Phi_{SBF_{+}^{-}}(\mathcal{X})\}
\end{equation*}
and
\begin{equation*}
\sigma_{LBW}(T)=\{\lambda\in\sigma_{a}(T): T-\lambda I\not\in\Phi_{SBF_{-}^{+}}(\mathcal{X})\},
\end{equation*}
respectively.
$T\in L(\mathcal{X})$ will be said B-Weyl, if $T$ is both upper and lower semi B-Weyl (equivalently $T$ is B-Fredholm operator of index zero).
The B-Weyl spectrum $\sigma_{BW}(T)$ of $T$ is defined by
\begin{equation*}
\sigma_{BW}(T)=\{\lambda\in\mathbb{C}: T-\lambda I\,\,\,\text{%
is not B-Weyl operator}\}.
\end{equation*}
Let $\Pi^{l}(T)$ denote the set of left pole of $T\in L(\mathcal{X})$.
\begin{equation*}
\Pi^{l}(T)=\{\lambda\in\sigma_{a}(T): asc(T-\lambda I)=d<\infty\,\,\,\text{and}\,\ R((T-\lambda I)^{d+1})\,\,\text{is closed}\}.
\end{equation*}
A strong version of the left polaroid property says that
$T\in L(\mathcal{X})$ is finitely left polaroid (respectively, finitely right polaroid)
if and only if every $\lambda\in iso\sigma_{a}(T)$
( respectively, $\lambda\in iso\sigma_{s}(T)$)
is left pole of $T$ and $\alpha(T-\lambda I)<\infty$ (respectively, right pole of $T$
and $\beta(T-\lambda I)<\infty$). Let $\Pi_{0}^{l}(T)$ (respectively, $\Pi_{0}^{r}(T)$)
denote the set of finite left poles (respectively, the set of finite right poles ) of $T.$
Then $T\in L(\mathcal{X})$ is finitely left polaroid (respectively, finitely right polaroid) if and only if
$iso \sigma_a(T)=\Pi_{0}^{l}(T)$ (respectively, $iso \sigma_a(T)=\Pi_{0}^{r}(T)$).\\
\indent For $T\in L(\mathcal{X})$ define
$$ \Delta (T)= \{ n \in \mathbb{N} : m \geq n, m\in  \mathbb{N} \Rightarrow R(T^n) \cap N(T) \subseteq R(T^m) \cap N(T) \}.$$
The degree of stable iteration is defined as $dis(T)=\inf \Delta (T) $ if $\Delta (T)\neq \emptyset,$ while
$dis(T)=\infty$ if $\Delta(T)=\emptyset.$ $T\in L(\mathcal{X})$
is said to be quasi-Fredholm of degree d, if there exists $d\in \mathbb{N} $ such that
$dis(T)=d,$ $R(T^n)$ is a closed subspace of $\mathcal{X}$ for each $n\geq d$ and $R(T)+N(T^n)$ is a closed subspace
of $\mathcal{X}.$ An operator $T\in L(\mathcal{X})$ is said to be semi-regular,
if $R(T)$ is closed and $N(T^{n})\subseteq R(T^{m})$ for all $m,n \in \mathbb{N} .$\\
\indent An important property in local spectral theory is the single valued extension property.
An operator $T\in L(\mathcal{X})$ is said to have the single valued extension property
at $\lambda_{0}\in\mathbb{C}$ (abbreviated SVEP at $\lambda_{0}$), if for every open
disc $\mathbb{D}$ centered at $\lambda_{0}$, the only analytic function $f:\mathbb{D}\rightarrow \mathcal{X}$
which satisfies the equation $(T-\lambda I)f(\lambda)=0$ for all $\lambda\in\mathbb{D}$
is the function $f\equiv0$. An operator $T\in L(\mathcal{X})$ is said to have SVEP if $T$ has SVEP at every $\lambda\in\mathbb{C}.$\\
Furthermore, for $T\in L(\mathcal{X})$ the quasi-nilpotent part of $T$ is defined by
\begin{equation*}
H_{0}(T)=\{x\in \mathcal{X}: \lim_{n\rightarrow\infty}\|T^{n}(\mathcal{X})\|^{\frac{1}{n}}=0\}.
\end{equation*}
It is easily seen that $N(T^{n})\subset H_{0}(T)$ for every $n\in\mathbb{N}$.
The analytic core of an operator $T\in L(\mathcal{X})$ is the subspace $K(T)$ defined as the set of all $x\in \mathcal{X}$ such that there exists a constant $c>0$ and a sequence of elements $x_{n}\in \mathcal{X}$ such that $x_{0}=x$, $Tx_{n}=x_{n-1},$ and $\|x_{n}\|\leq c^{n}\|x\|$ for all $n\in\mathbb{N}$, the spaces $K(T)$ are hyperinvariant under $T$ and satisfy $K(T)\subset R(T^{n})$, for every $n\in\mathbb{N}$ and  $T(K(T))=K(T),$  see \cite{P1} for information on $H_{0}(T)$ and $K(T)$.

\section{Left polaroid generalized derivation}

We begin this section by recalling some results concerning spectra of generalized derivations.\\
Let $\mathcal{X}$ and $\mathcal{Y}$ be two Banach spaces and consider
$A\in L(\mathcal{X})$ and $B\in L(\mathcal{Y}).$ Let
$\delta_{A,B} \in L(L(\mathcal{Y},\mathcal{X}))$ the generalized derivation induced by $A$
and $B,$ i.e., $$\delta_{A,B}(X)=(L_{A}-R_{B})(X)=AX-XB \mbox{ where } X\in L(\mathcal{Y},\mathcal{X}).$$
According to \cite[Theorem 3.5.1]{LN}, we have that $$\sigma_{a}(\delta_{A,B})=\sigma_{a}(A)-\sigma_{s}(B).$$
and it is not difficult to conclude that
$$iso \sigma_a (\delta_{A,B})= (iso\sigma_a (A)- iso\sigma_a (B^*))\setminus acc\sigma_a(\delta_{A,B}).$$

The following results concerning upper semi Fredholm spectrum and Browder essential approximate point spectrum
of generalized derivation was proved in \cite{BA, LZ}. It will be used in the sequel.

\begin{lemma} \label{propostion1}
Let $\mathcal{X}$ and $\mathcal{Y}$ be two Banach spaces and consider $A\in L(\mathcal{X})$
and $B\in L(\mathcal{Y}),$ then the following statements hold.
\begin{enumerate}
  \item[i)] $\sigma_{SF_{+}}(\delta_{A,B})=(\sigma_{SF_{+}}(A)-\sigma_{s}(B))\cup(\sigma_{a}(A)-\sigma_{SF_{-}}(B)).$
  \item[ii)] $\sigma_{ab}(\delta_{A,B})=(\sigma_{ab}(A)-\sigma_{s}(B))\cup(\sigma_{a}(A)-\sigma_{ab}(B^{*})).$
  \end{enumerate}
\end{lemma}

The following lemma  concerning the Weyl essential approximate point spectrum of generalized derivation
will be used in the sequel.

\begin{lemma}\label{lema2}
Let $\mathcal{X}$ and $\mathcal{Y}$ be two Banach spaces and consider $A\in L(\mathcal{X})$
and $B\in L(\mathcal{Y}),$ then
 $$\sigma_{aw}(\delta_{A,B})\subseteq(\sigma_{aw}(A)-\sigma_{s}(B))\cup(\sigma_{a}(A)-\sigma_{aw}(B^{*})).$$
\end{lemma}
\begin{proof}
Let $\lambda\not\in(\sigma_{aw}(A)-\sigma_{s}(B))\cup(\sigma_{a}(A)-\sigma_{aw}(B^{*})).$
If $\mu_{i}\in\sigma_{a}(A)$ and $\nu_{i}\in\sigma_{s}(B)$ are such that $\lambda=\mu_{i}-\nu_{i},$
then $\mu_{i}\not\in\sigma_{SF_{+}}(A)$ and $\nu_{i}\not\in\sigma_{SF_{-}}(B),$ hence from
statement i) of Lemma \ref{propostion1} $\lambda\not\in\sigma_{SF_{+}}(\delta_{A,B})$.
Now, we will prove that $$ind(\delta_{A,B}-\lambda I)\leq0.$$
Suppose to the contrary that $ind(\delta_{A,B}-\lambda I)>0$, then $\lambda\not\in\sigma_{e}(\delta_{A,B}).$ It follows from \cite[Corollary 3.4]{E}
that $$\lambda=\mu_{i}-\nu_{i}\,\,\,\ (1\leq i\leq n),$$ where
$\mu_{i}\in iso\sigma(A)$ for $1\leq i\leq m$ and $\nu_{i}\in iso\sigma(B),$ for
$m+1\leq i\leq n$. We have that
$ind(\delta_{A,B}-\lambda I)$ is equal to
$$\sum_{j=m+1}^{n}dim H_{0}(B-\nu_{j})ind(A-\mu_{j})-\sum_{k=1}^{m}dim H_{0}(A-\mu_{k})ind(B-\nu_{k}).$$
Since $\mu_{i}\in iso\sigma(A),$ for $1\leq i\leq m$ and
$\nu_{i}\in iso\sigma(B),$ for $m+1\leq i\leq n$, it follows that
$dim H_{0}(A-\mu_{j})$ is finite, for $1\leq j\leq m$ and $dim H_{0}(B-\nu_{k})$ is finite,
for $m+1\leq k\leq n$ and we have also $ind(A-\mu_{i})\leq0$ and $ind(B-\nu_{j})\geq0.$
Thus $ind(\delta_{A,B}-\lambda I)\leq0.$ This a contradiction.
Hence $\lambda\not\in\sigma_{aw}(\delta_{A,B}).$
\end{proof}

\indent According to \cite{BD} left polaroid operator (respectively, right polaroid ) operator satisfies property
$(\mathcal{P}_{l}),$ (respectively, $(\mathcal{P}_{r})$), if it is left polar at every
$\lambda \in iso\sigma_a(T)$ (respectively, right polar at every $\lambda \in iso\sigma_s(T)$
which satisfies property $(\mathcal{P}_{l}),$(respectively property $(\mathcal{P}_{r}).$\\
\indent The Following Lemma is the dual version of \cite[Lemma 3.1]{BD}.

\begin{lemma}\label{lemmaD}
Let $\mathcal{X}$ a Banach space.
If $T\in  L(\mathcal{X})$ is right polaroid and satisfies property $(\mathcal{P}_{r}),$
then for every $\lambda \in iso\sigma_s(T)$ there exists T-invariant closed
subspaces $N_1$ and $N_2$ such that
$\mathcal{X}=N_1 \oplus N_2,$ $(T-\lambda)|_{N_1}$ is nilpotent of order $d(\lambda)$ and
$(T-\lambda I)|_{N_2}$ is surjective, where $d(\lambda)$ is the order of the right pole at $\lambda.$
Moreover, $K(T-\lambda I)=R((T- \lambda I)^{d(\lambda)}).$
\end{lemma}
\begin{proof}
From the hypothesis $T-\lambda$ is quasi-Fredholm of degree $d(\lambda)$
and the closed subspaces $N((T-\lambda I)^{d(\lambda)})+R(T-\lambda)$ is complemented in
$\mathcal {X}.$ Since $T\in  L(\mathcal{X})$ is right polaroid and satisfies property $(\mathcal{P}_{r}),$
then $N(T-\lambda) \cap R((T-\lambda)^{d(\lambda)})$ is complemented in $\mathcal {X}.$
From \cite[ Theorem 5]{Mu}, there exists T-invariant closed subspaces $N_1$ and $N_2$
such that $\mathcal{X}=N_1 \oplus N_2,$ $(T-\lambda)|_{N_1}$ is nilpotent of order $d(\lambda)$ and
$(T-\lambda I)|_{N_2}$ is semi-regular.
Since $dsc(T-\lambda I)=d(\lambda),$ then
the semi-regular operator $(T-\lambda I)|_{N_2}$ is surjective. Since
$K(T-\lambda I)=K((T-\lambda I)|{N_1}) \oplus K((T-\lambda I)|{N_2})= 0 \oplus N_2=N_2,$ then
we can conclude from \cite[Theorem 2.7]{P2} that $K(T-\lambda I)=R((T- \lambda I)^{d(\lambda)}.$
\end{proof}
Next follows the main result of this section.
\begin{theorem}\label{theorem1}
Let $\mathcal{X}$ and $\mathcal{Y}$ be two Banach spaces and let $A\in L(\mathcal{X})$
be left polaroid and $B\in L(\mathcal{Y})$ be right polaroid.
If $A$ satisfies property $(\mathcal{P}_l)$ and $B$ satisfy property $(\mathcal{P}_r),$ then
\begin{center}$\delta_{A,B}$ is left polaroid.\end{center}
\end{theorem}
\begin{proof}
Let $\lambda\in iso\sigma_{a}(\delta_{A,B}),$ then
there exist $\mu\in\sigma_{a}(A)$ and $\nu\in\sigma_{s}(B)$ such that $\lambda=\mu-\nu$,
it follows from statement that $\mu\in iso\sigma_{a}(A)$ and $\nu\in iso\sigma_{s}(B)=iso\sigma_{a}(B^{*}).$
Since $A$ is left polaroid, then there exist $A$-invariant closed subspaces
$M_{1}$ and $M_{2}$ such that $\mathcal{X}=M_1 \oplus M_2,$
$(A- \mu I)|_{M_{1}}= A_{1}- \mu I|_{M_1}$ is nilpotent of order $d_1$ where
where $d_1=d(\mu)$ is the order of the left pole of $A$ at $\mu$
and that $(A- \mu I)|_{M_2}=A_{2}- \mu I|_{M_2}$ is bounded below.
Also, since $B$ is right polaroid,
then there exists $B$-invariant closed subspaces $N_1$ and $N_2$ such that
$\mathcal{Y}=N_1 \oplus N_2,$ $(B-\nu)|_{N_1}=B_{1}-\nu I |_{N_1}$ is nilpotent of order $d_2$
where $d_2=d(\nu)$ is the order of the right pole of $B$ at $\nu$ and $(B-\nu I)|_{N_2}=B_{2}-\nu I|_{N_2}$ is surjective.
Let $d=d_{1}+d_{2}$
and $X\in L(N_{1}\oplus N_{2},M_{1}\oplus M_{2})$ have the representation $X=[X_{kl}]_{k,l=1}^{2}$.
We will prove that $asc(\delta_{A,B}-\lambda I)$ is finite.
We have
\begin{eqnarray*}
(\delta_{A_1,B_1}-\lambda I)^{d+1}&=&(L_{A_1-\mu I}-R_{B_1-\nu I})^{d+1} \\
&=&\sum_{k=0}^{d+1}(-1)^{d+1-k}\left(
                    \begin{array}{c}
                      d+1 \\
                      k\\
                    \end{array}
                  \right)(L_{A_1-\mu I})^{k}(R_{B_1-\nu I})^{d+1-k}=0,
\end{eqnarray*} that is
$\delta_{A_{1},B_{1}}-\lambda I$ is nilpotent of order $d,$ and hence $asc(\delta_{A_{1},B_{1}}-\lambda I)<\infty .$\\
On the other hand,
 \begin{eqnarray*}
 (\delta_{A_{1},B_{2}}-\lambda I)^{d+1}(X_{12})&=&\sum_{k=0}^{d_{1}-1}(-1)^{d_{1}-1-k}\left(
                    \begin{array}{c}
                      d_{1}-1 \\
                      k\\
                    \end{array}
                  \right)(L_{A_{1}-\mu I})^{k}(R_{B_{2}-\nu I})^{d_{1}-1-k}(X_{12})\\
&=&(\delta_{A_{1},B_{2}}-\lambda I)^{d}(X_{12}),
\end{eqnarray*}similarly we get
$(\delta_{A_{2},B_{1}}-\lambda I)^{d+1}(X_{21})=(\delta_{A_{2},B_{1}}-\lambda I)^{d}(X_{21}).$ Thus
$$asc(\delta_{A_{1},B_{2}}-\lambda I)<\infty \mbox{ and } asc(\delta_{A_{2},B_{1}}-\lambda I)<\infty .$$
Now, we will prove that $0\not\in\sigma_{a}(\delta_{A_{2}-\mu I|_{M_2},B_{2}-\nu I|_{N_2}}).$
For this, it suffices to prove that $\sigma_{a}(A_{2}-\mu I|_{M_2}  ) \cap \sigma_{s}(B_{2}-\nu I|_{N_2})= \emptyset.$
Suppose that there exists a complex number $\alpha $ such that
$\alpha \in  \sigma_{a}(A_{2}-\mu I|_{M_2})  \cap \sigma_{s}(B_{2}-\nu I|_{N_2}),$ then
$\alpha \in\sigma_{a}(A_{2}-\mu I|_{M_2})$ and $ \alpha \in\sigma_{s}(B_{2}-\nu I|_{N_2}),$
from \cite[Theorem 2.48]{P1}, $0\in\sigma_{a}(A_{2}-(\mu +\alpha)I|_{M_2}  )$ and $0 \in\sigma_{s}(B_{2}-(\nu + \alpha)I|_{N_2}).$
Since $(\mu +\alpha)$ is isolated in the approximate point spectrum of $A$ and $(\nu +\alpha)$ is isolated in the surjective spectrum of $B,$
then by the hypothesis $A$ is left polaroid which satisfies property $(\mathcal{P}_l)$
and $B$ is right polaroid which satisfies property $(\mathcal{P}_r),$ we conclude that
$$(A- (\mu +\alpha )I)|_{M_2}=A_{2}- (\mu + \alpha)I|_{M_2}$$ is bounded below and
$$(B-(\nu+\alpha)I)|_{N_2}=B_{2}-(\nu +\alpha)I|_{N_2}$$ is surjective. That is
$$0 \notin \sigma_{a}(A_{2}-(\mu +\alpha)I|_{M_2}) \mbox{ and } 0 \notin \sigma_{s}(B_{2}-(\nu + \alpha)I|_{N_2}).$$
This is a contradiction, hence $0\not\in\sigma_{a}(\delta_{A_{2}-\mu I|_{M_2},B_{2}-\nu I|_{N_2}}).$
Finally, $$asc(\delta_{A_2,B_2}-\lambda I)\leq d<\infty,$$ and hence, $$asc(\delta_{A,B}-\lambda I)\leq d<\infty.$$
Now,we prove that $(\delta_{A,B}-\lambda I)^{d+1}(L(\mathcal{Y},\mathcal{X}))$ is closed. Since
$0\not\in\sigma_{a}(\delta_{A_{2},B_{2}}-\lambda I),$ then from \cite[Lemma 1.1]{AA}
$(\delta_{A_{2},B_{2}}-\lambda I)^{d+1}(L(N_{2},M_{2})$ is closed.
We have that $\delta_{A_{1},B_{1}}-\lambda I$ is nilpotent of order $d,$ then by \cite[Theorem 2.7]{OBO}
it follows that $$(\delta_{A_{1},B_{1}}-\lambda I)^{d+1}(L(N_{1},M_{1})) \mbox{ is closed }.$$
From the fact that $0\not\in\sigma_{a}(\delta_{A_{i},B_{j}}-\lambda I)$ and \cite[Lemma 1.1]{AA},
it follows that $$(\delta_{A_{i},B_{j}}-\lambda I)^{d+1}(L(N_{j},M_{i}) \mbox{
is closed for }1\leq i,j\leq2 \mbox{ and }i\neq j,$$ consequently
$(\delta_{A,B}-\lambda I)^{d+1}(L(\mathcal{X},\mathcal{Y}))$ is closed. Hence
$\lambda$ is a left pole of $\delta_{A,B}$ which means that $\delta_{A,B}$ is left polaroid.
\end{proof}
In the case of Hilbert spaces, we have the following corollary.
\begin{corollary}
Let $H$ and $K$ be Hilbert spaces and let $A\in L(H)$ and $B \in L(K).$
If $A$ and $B^*$ are left polaroid, then
\begin{center}$\delta_{A,B}$ is left polaroid.\end{center}
\end{corollary}

In the following Theorem, we characterize finitely left polaroid generalized derivation.

\begin{theorem}
Let $\mathcal{X}$ and $\mathcal{Y}$ be two Banach spaces and let $A\in L(\mathcal{X})$
and $B\in L(\mathcal{Y})$. Then
$A$ and $B^*$ are finitely left polaroid operators if and only if
\begin{center}$\delta_{A,B}$ is finitely left polaroid.\end{center}
\end{theorem}
\begin{proof}
Let $\lambda\in iso\sigma_{a}(\delta_{A,B}),$ then there exist
$\mu\in\sigma_{a}(A)$ and $\nu\in\sigma_{s}(B)$ such that
$\lambda=\mu-\nu,$ hence we have
$\mu\in iso\sigma_{a}(A)$ and $\nu\in iso\sigma_{s}(B)=iso\sigma_{a}(B^{*}).$
Suppose that  $A$ and $B^{*}$ are finitely left polaroid, then from \cite[Corollary 2.2]{R1} we have
that $\mu\not\in\sigma_{ab}(A)$ and $\nu\not\in\sigma_{ab}(B^{*})$, applying statement ii) of Lemma \ref{propostion1},
we get $\lambda\not\in\sigma_{ab}(\delta_{A,B}),$
hence by \cite[Corollary 2.2]{R1} $\delta_{A,B}$ is finitely left polaroid.
Conversely, suppose that $\delta_{A,B}$ is finitely left polaroid and prove that
$A$ and $B^{*}$ are finitely left polaroid.
For this, let $\mu\in iso\sigma_{a}(A)$ and $\nu\in iso\sigma_{a}(B^{*}),$ then
$\lambda=\mu-\nu \in iso\sigma_{a}(\delta_{A,B}).$ Since
$\delta_{A,B}$ is finitely left polaroid, then by \cite[Corollary 2.2]{R1}
$\lambda=\mu - \nu \not\in\sigma_{ab}(\delta_{A,B}),$ and hence by statement ii) of Lemma \ref{propostion1}
$\mu\notin\sigma_{ab}(A)$ and $\nu\notin\sigma_{ab}(B^{*}).$ We conclude from
\cite[Corollary 2.2]{R1} that $A$ and $B^{*}$ are finitely left polaroid.
\end{proof}

\section{Consequences on Weyl's type theorem}

For  $T\in L(\mathcal{X})$, let $E^{a}(T)=\{\lambda\in iso\sigma_{a}(T): 0<\alpha(T-\lambda I)\}$ and
$E_{0}^{a}(T)=\{\lambda\in E^{a}(T): \alpha(T-\lambda I)<\infty\}.$
Recall that $T$ is said to satisfy a-Browder's theorem (respectively, generalized a-Browder's theorem)
if $\sigma_{a}(T)\backslash\sigma_{aw}(T)=\Pi_{0}^{l}(T)$ (respectively, $\sigma_{a}(T)\backslash\sigma_{UBW}(T)=\Pi^{l}(T)$).
From \cite[Theorem 2.2]{AZ6} we have $T$ satisfies a-Browder's theorem if and only if
$T$ satisfies generalized a-Browder's theorem.
$T$ is said to satisfy a-Weyl's theorem (respectively, generalized a-Weyl's theorem) if
$\sigma_{a}(T)\backslash\sigma_{aw}(T)=E_{0}^{a}(T)$
(respectively, $\sigma_{a}(T)\backslash\sigma_{UBW}(T)=E^{a}(T)$).\\
For  $T\in L(\mathcal{X})$, let $E(T)=\{\lambda\in iso\sigma(T): 0<\alpha(T-\lambda I)\}$ and
$E_{0}(T)=\{\lambda\in E(T): \alpha(T-\lambda I)<\infty\}.$ Recall that
$T$ is said to satisfy Weyl's theorem (respectively, generalized Weyl's theorem) if
$\sigma(T)\backslash\sigma_{W}(T)=E_{0}(T)$ (respectively, $\sigma(T)\backslash\sigma_{BW}(T)=E(T)$).
We know that $T$ satisfies generalized a-Weyl's theorem implies that
$T$ satisfies a-Weyl's theorem and this implies that $T$ satisfies Weyl's theorem.
Next generalized a-Weyl's theorem for $\delta_{A,B}$ will be studied.

\begin{theorem}\label{thmanswer}
Let $\mathcal{X}$ and $\mathcal{Y}$ be two Banach spaces and let $A\in L(\mathcal{X})$
and $B\in L(\mathcal{Y})$. Suppose that $A$ and $B^{*}$ satisfy a-Browder's theorem.
If $A$ is left polaroid and satisfies property $(\mathcal{P}_l)$ and $B$ is right polaroid
and satisfies $(\mathcal{P}_r)$, then the following assertions are equivalent.
\begin{itemize}
\item[i)] $\delta_{A,B}$ satisfies generalized a-Weyl's theorem.
\item[ii)] $\sigma_{aw}(\delta_{A,B})=(\sigma_{aw}(A)-\sigma_{s}(B))\cup(\sigma_{a}(A)-\sigma_{aw}(B^{*})).$
\end{itemize}
\end{theorem}
\begin{proof} If $A$ and $B^{*}$ satisfy a-Browder theorem, then they satisfy generalized a-Browder theorem,
by \cite[Theorem 4.2]{BA} it follows that $\delta_{A,B}$ satisfies generalized a-Browder's theorem if and only if $\sigma_{aw}(\delta_{A,B})=(\sigma_{aw}(A)-\sigma_{s}(B))\cup(\sigma_{a}(A)-\sigma_{aw}(B^{*})).$ That is $\sigma_{a}(\delta_{A,B})\backslash\sigma_{UBW}(\delta_{A,B})=\Pi^{l}(\delta_{A,B}).$
Since  $A$ is left polaroid and $B$ is right polaroid,
then from Theorem \ref{theorem1} $\delta_{A,B}$ is left polaroid,
consequently $\Pi^{l}(\delta_{A,B})=E^{a}(\delta_{A,B}).$
Thus $\delta_{A,B}$ satisfies generalized a-Weyl's theorem.
The reverse implication is obvious from the implication $\delta_{A,B}$ satisfies generalized a-Weyl's theorem implies
$\delta_{A,B}$ satisfies generalized a-Browder's theorem
\end{proof}
In the case of Hilbert spaces operators, we have the following corollaries.
\begin{corollary}
Let $H$ and $K$ be two Hilbert spaces and let $A\in L(H)$
and $B\in L(K)$. Suppose that $A$ and $B^{*}$ satisfy a-Browder's theorem.
If $A$ is left polaroid and $B$ is right polaroid, then the following assertions are equivalent.
\begin{itemize}
\item[i)] $\delta_{A,B}$ satisfies generalized a-Weyl's theorem.
\item[ii)] $\sigma_{aw}(\delta_{A,B})=(\sigma_{aw}(A)-\sigma_{s}(B))\cup(\sigma_{a}(A)-\overline{\sigma_{aw}(B^{*})}).$
\end{itemize}
\end{corollary}

\begin{corollary}\label{coro}
Let $\mathcal{X}$ and $\mathcal{Y}$ be two Banach spaces and let $A\in L(\mathcal{X})$
and $B\in L(\mathcal{Y})$. Suppose that $A$ and $B^{*}$ satisfy a-Browder's theorem.
If $A$ is left polaroid and satisfies property $(\mathcal{P}_l)$ and
$B$ is right polaroid and satisfies property $(\mathcal{P}_r)$, then the following assertions are equivalent.
\begin{itemize}
\item[i)] $\delta_{A,B}$ has SVEP at $\lambda\not\in\sigma_{UBW}(\delta_{A,B}).$
\item[ii)] $\delta_{A,B}$ satisfies a-Browder's theorem.
\item[iii)] $\delta_{A,B}$ satisfies a-Weyl's theorem.
\item[iv)] $\delta_{A,B}$ satisfies generalized a-Weyl's theorem.
\item[v)] $\sigma_{aw}(\delta_{A,B})=(\sigma_{aw}(A)-\sigma_{s}(B))\cup(\sigma_{a}(A)-\sigma_{aw}(B^{*})).$
\end{itemize}
\end{corollary}
\begin{proof}
$(i)\Leftrightarrow(ii)$ follows from \cite[Theorem 2.1]{AZ8}, $(iii)\Leftrightarrow(iv)$ follows from \cite[Theorem 3.7]{AA} and $(iv)\Leftrightarrow(v)$ follows from Theorem \ref{thmanswer}.
\end{proof}
In the following result, we give sufficient conditions for $\delta_{A,B}$ to satisfy
a-Browder's theorem.
\begin{theorem}\label{theo}
Let $\mathcal{X}$ and $\mathcal{Y}$ be two Banach spaces and let $A\in L(\mathcal{X})$
and $B\in L(\mathcal{Y})$.
If $A$ has SVEP on the complement of $\sigma_{SF_{+}}(A)$ and
$B$ has SVEP on the complement of $\sigma_{SF_{-}}(B)$, then
\begin{center}
$\delta_{A,B}$ satisfies a-Browder's theorem.
\end{center}
\end{theorem}
\begin{proof}
Let $\lambda\in\sigma_{a}(\delta_{A,B})\backslash\sigma_{aw}(\delta_{A,B}),$ then $\lambda\in\sigma_{a}(\delta_{A,B})\backslash\sigma_{SF_{+}}(\delta_{A,B}),$ from statement i) of Lemma \ref{propostion1} there exist $\mu\in\sigma_{a}(A)\backslash\sigma_{SF_{+}}(A)$ and $\nu\in\sigma_{s}(B)\backslash\sigma_{SF_{-}}(B)$ such that $\lambda=\mu-\nu$. Since $A$ has SVEP at $\mu\not\in\sigma_{SF_{+}}(A)$ and $B$ has SVEP at $\nu\not\in\sigma_{SF_{-}}(B)$, it follows from \cite[Corollary 2.2]{R1} that $\mu\not\in\sigma_{ab}(A)$ and $\nu\not\in\sigma_{ab}(B^{*})$, applying
statement ii) of Lemma \ref{propostion1} we get $\lambda\not\in\sigma_{ab}(\delta_{A,B})$ this is equivalent to $\lambda\in\Pi_{0}^{l}(\delta_{A,B}).$ Hence $\delta_{A,B}$ satisfy a-Browder's theorem.
\end{proof}

\section{Application}
A Banach space operator $T\in L(\mathcal{X})$ is said to be hereditary normaloid, $T\in\mathcal{HN}$, if every part of $T$ (i.e., the restriction of $T$ to each of its invariant subspaces) is normaloid (i.e., $\|T\|$ equals the spectral radius $r(T)$), $T\in\mathcal{HN}$ is totally hereditarily normaloid $\mathcal{THN}$ if also the inverse of every invertible part of $T$ is normaloid and $T$ is completely totally hereditarily normaloid
(abbreviated $T\in\mathcal{CHN}$), if either $T\in\mathcal{THN}$ or $T-\lambda I\in\mathcal{HN}$ for every complex number $\lambda.$ The class $\mathcal{CHN}$ is large. In particular, let $H$ a Hilbert space and $T\in L(H)$ a Hilbert space operator. If $T$ is hyponormal ($T^{*}T\geq TT^{*}$) or  p-hyponormal ($(T^{*}T)^{p})\geq(TT^{*})^{p}$) for some ($0<p\leq1$) or w-hyponormal  $((|T^{*}|^{\frac{1}{2}}|T||T^{*}|^{\frac{1}{2}})^{\frac{1}{2}}\geq|T^{*}|),$ then $T$ is in $\mathcal{THN}.$ Again totaly *-paranormal operators ($\|(T-\lambda I)^{*}x\|^{2}\leq\|(T-\lambda I)x\|^{2}$ for every unit vector $x$) are $\mathcal{HN}$ operator and paranormal operators ($\|Tx\|^{2}\leq\|T^{2}x\|\|x\|,$ for all unit vector $x$) are $\mathcal{THN}$ operators.
It is proved in \cite{D2} that if $A,B^{*}\in L(H)$ are hyponormal, then generalized Weyl's theorem holds for $f(\delta_{A,B})$ for every $f\in\mathcal{H}(\sigma(\delta_{A,B})),$ where $\mathcal{H}(\sigma(\delta_{A,B}))$ is the set of all analytic functions defined on a neighborhood of $\sigma(\delta_{A,B}),$ this result was extended to log-hyponormal or p-hyponormal operators in \cite{Du8} and \cite{LB}. Also in \cite{CH} and \cite{LBA} it is shown that if $A,B^{*}\in L(H)$ are w-hyponormal operators, then Weyl's theorem holds for $f(\delta_{A,B})$ for every $f\in\mathcal{H}(\sigma(\delta_{A,B})).$
Let $\mathcal{H}_{c}(\sigma(T))$ denote the space of all analytic functions defined on a neighborhood of $\sigma(T)$ which is non constant on each of the components of its domain. In the next results we can give more.
\begin{theorem}\label{theorem4}
Suppose that $A,B\in L(H)$ are $\mathcal{CHN}$ operators, then \begin{center}
$\delta_{A,B}$ satisfies a-Browder's theorem.
\end{center}
\end{theorem}
\begin{proof}
Since $A$ and $B$ are $\mathcal{CHN}$ operators, it follows from \cite[Corollary 2.10]{D5} that $A$ has SVEP on the complement of $\sigma_{SF_{+}}(A)$ and $B$ has SVEP on the complement of $\sigma_{SF_{-}}(B),$ then by Theorem \ref{theo} a-Browder's theorem holds for $\delta_{A,B}$.
\end{proof}
\begin{corollary}\label{corolla} Suppose that $A,B\in L(H)$ are $\mathcal{CHN}$ operators, then the following assertions are equivalent.
\begin{itemize}
\item[i)] $\delta_{A,B}$ has SVEP at $\lambda\not\in\sigma_{UBW}(\delta_{A,B})$
\item[ii)] $\delta_{A,B}$ satisfies a-Browder's theorem.
\item[iii)] $\delta_{A,B}$ satisfies a-Weyl's theorem.
\item[iv)] $\delta_{A,B}$ satisfies generalized a-Weyl's theorem.
\item[v)] $\sigma_{aw}(\delta_{A,B})=(\sigma_{aw}(A)-\sigma_{s}(B))\cup(\sigma_{a}(A)-\overline{\sigma_{aw}(B^{*})}).$
\end{itemize}
\end{corollary}

\begin{proof}
Since $A$ and $B$ are $\mathcal{CHN}$ operators, it follows from \cite[Corollary 2.15]{D5} that $A$, $B$, $A^{*}$ and $B^{*}$ satisfy a-Browder's theorem.
By \cite[Proposition 2.1]{D5}, we conclude that $A$ and $B^{*}$ are left polaroid. Now,
then the equivalences follows from Corollary \ref{coro}.
\end{proof}
\begin{corollary} Suppose that $A,B\in L(H)$ are $\mathcal{CHN}$ operators, then \begin{center}
$f(\delta_{A,B})$ satisfies generalized a-Browder's theorem, for every $f\in\mathcal{H}_{c}(\sigma(\delta_{A,B}))$.\end{center}
\end{corollary}
\begin{proof}
By Corollary \ref{corolla} and \cite[Corollary 3.5]{D7}, we get generalized a-Browder's theorem holds for $f(\delta_{A,B})$.
\end{proof}
\begin{corollary} Suppose that $A,B\in L(H)$ are $\mathcal{CHN}$ operators, then \begin{center}
$f(\delta_{A,B})$ satisfies generalized a-Weyl's theorem, for every $f\in\mathcal{H}_{c}(\sigma(\delta_{A,B}))$.\end{center}
\end{corollary}
\begin{proof}
By \cite[Proposition 2.1]{D5} and Theorem \ref{theorem1},
 we get $\delta_{A,B}$ is left polaroid and from Corollary \ref{corolla} we have
 $\delta_{A,B}$ satisfies generalized a-Weyl's theorem,
 apply \cite[Theorem 3.14]{D7} we get generalized a-Weyl's's theorem holds for $f(\delta_{A,B})$.
\end{proof}

\end{document}